\documentclass[a4paper]{amsart}  

\usepackage{amssymb} \usepackage{amscd} \usepackage{times}

\def\zbb{\mathbb{Z}}  
  
  \def\phi{\varphi}
 \def\p1{{\mathbb{P}^1_\zbb}}

\begin{document}

\title{Effet d'une perturbation non lin\'eaire sur l'obtention d'une estimation uniforme. }
\author{Samy Skander Bahoura}

\address{6, rue Ferdinand Flocon, 75018 Paris, France. }
              
\email{samybahoura@yahoo.fr, bahoura@ccr.jussieu.fr} 

\date{}

\maketitle

\begin{abstract}

We consider the equation $ \Delta u= Vu^{(n+2)/(n-2)}+Wu^{n/(n-2)} $ and we give some minimal conditions on $ \nabla V $ and $ \nabla W $ to have an uniform estimate for their solutions.
 
If we replace $ Wu^{n/(n-2)} $ by $ W u $ in the previous equation, we have an uniform estimate for radial solutions.

\end{abstract}

\bigskip

\bigskip

\begin{center}  1. INTRODUCTION ET R\'ESULTATS.
\end{center}

\bigskip

Nous notons $ \Delta=-\sum \partial_{ii} $ le laplcien g\'eom\'etrique sur $ {\mathbb R}^n, n\geq 3 $.

\bigskip

Consid\'erons sur un ouvert $ \Omega $ de $ {\mathbb R}^n, n\geq 3 $, l'\'equation suivante:

$$ \Delta u =V u^{(n+2)/(n-2)}+ Wu^{n/(n-2)} \qquad (E) $$

o\`u $ V $ et $ W $ sont deux fonctions r\'eguli\`eres.

\bigskip

On suppose que:

$$ 0 < a \leq V(x) \leq b, \,\, ||\nabla V||_{L^{\infty}} \leq A \qquad (C_1) $$

$$ 0 < c \leq W(x) \leq d, \,\, ||\nabla V||_{L^{\infty}} \leq B \qquad (C_2) $$

{\bf Probl\`eme:} Quelle conditions minimales peut on imposer \`a $ \nabla V $ et $ \nabla W $ pour avoir une esimation uniforme du type $ \sup \times \inf $ pour les solutions de l'\'equation $ (E) $ ?

\bigskip

Notons que lorsque $ W \equiv 0 $, l'\'equation $ (E) $ est la c\'el\`ebre \'equation de la courbure scalaire prescrite sur un ouvert de l'\'espace euclidien de dimension $ n \geq 3 $.

\bigskip

Dans ce cas, il existe beaucoup de r\'esultats concernant les solutions de cette \'equation, voir par exemple, [B], [C-L 1].

\bigskip

Lorsque $ \Omega = {\mathbb S}_n $ avec l'\'equation correspondante (courbure scalaire),
 YY. Li donne des conditions de platitude suffisante pour avoir une majoration de l'\'energie ainsi que l'existence de points dit isol\'es simples, voir [L1], [L2].

\bigskip

Dans [C-L 2]  Chen et Lin mettent en \'evidence un contre example confirmant l'importance des hypoth\`eses de Li.

\bigskip

Notons que dans [C-L 1], il existe des r\'esultats concernant les in\'egalit\'es de Harnack du type $ \sup \times \inf $ avec des conditions de platitudes similaires \`a celles de Li pour une equation du type:

$$ \Delta u= Vu^{(n+2)/(n-2)}+ g(u) $$

avec $ g $ une fonction r\'eguli\`ere \'equivalente $ t^{\alpha} $, $ 1 \leq \alpha < \dfrac{n+2}{n-2} $.

\bigskip

Notons que dans ce travail, aucune borne a priori sur l'\'energie n'est impos\'ee. On utilise la technique blow-up et de d\'eplacement de plan dite " Moving-Plane" invent\'ee par Alexandrov et d\'evelopp\'ee par Gidas-Ni-Nirenberg, voir [ G-N-N].

\bigskip

Notons que la m\'ethode "moving-Plane" est souvent utilis\'ee pour d\'eterminer si les solutions d'une EDP sont sym\'etriques ou dans la recherche de forme explicite de certaines solutions d'\'equations aux d\'eriv\'ees partielles.

\bigskip

Ici, nous avons:

\bigskip

{ \it {\bf Th\'eor\`eme 1}. Pour tout $ a,b,c,d >0 $, pour toutes suites $ (A_i), (B_i) $ telles que $ A_i \to 0 $ et $ B_i \to 0 $ et pour tout compact $ K $ de $ \Omega $, il existe une constante positive $ c=c[a,b,c,d, (A_i), (B_i), K, \Omega,n] $ telle que:

$$ \sup_K u_i \times \inf_{\Omega} u_i \leq c, \,\,\,\, {\rm ( pour \,\, i \,\, assez \,\, grand)}  $$

pour toute suite $ (u_i)_i $ solutions de $ (E) $ relativement \`a $ (V_i) $ et $ (W_i) $ v\'erifiant les conditions $ (C_1) $ et $ (C_2) $.}

\bigskip

On se place sur la boule unit\'e de $ {\mathbb R}^n $ ( $\Omega=B_1(0) $) et on s'occupe de l'\'equation suivante:

$$ \Delta u_i=V_i{u_i}^{(n+2)/(n-2)}+ W_i u_i \qquad (E') $$

On suppose que $ u_i $ et $ V_i $ sont radiales:

$$ 0 < a \leq V_i(r) \leq b\,\,\,{\rm et} \,\,\,|V_i(r)-V_i(r')| \leq A_i|r^2-r'^2| \,\,\, {\rm avec} \,\,\, A_i \to 0 \qquad (C_3) $$

$$ 0 < c \leq W_i(r) \leq d \,\,\,{\rm et} \,\,\, |W_i'(r)|\leq B_i \,\,\, {\rm avec} \,\,\, B_i \to 0 \qquad (C_4) $$

Nous avons:

\bigskip

{\it {\bf Th\'eor\`eme 2}. Pour tout $ a,b,c,d >0 $, pour toutes suites $ (A_i) $ et $ (B_i) $, il existe une constante positive $ c=c[a,b,c,d, (A_i), (B_i),n] $ telle que:

$$ u_i(0) \times u_i(1) \leq c \,\,\,\, {\rm  (pour \,\, i \,\, assez \,\, grand)} $$

pour toute suite $ (u_i) $ solution de $ (E') $ relativement \`a $ (V_i) $ et $ (W_i) $ v\'erifiant $ (C_3) $ et $ (C_4) $.}

\bigskip

\begin{center}  2. PREUVES DES TH\'EOR\`EMES.
\end{center}

\bigskip

\underbar {\bf Preuve du Th\'eor\`eme 1}

\bigskip

Soit $ x_0 $ un point de $ \Omega $ et $ (u_i)_i $ une suite de fonctions sur $ \Omega $ telles que,

$$ \Delta u_i = V_i {u_i}^{(n+2)/(n-2)}+W_i {u_i}^{n/(n-2)}, \,\, u_i>0 $$

On raisonne par l'absurde, en supposant que $ \sup \times \inf $ n'est pas born\'e.

\bigskip

On suppose que:

$ \forall \,\, c,R >0 \,\, \exists \,\, (u_{c,R,j})_j $ solution de $ (E) $ telle que:

$$ R^{n-2} \sup_{B(x_0,R)} u_{c,R,j}\times \inf_{ \Omega } u_{c,R,j} \geq c, \qquad (H) $$

\underbar {\bf Proposition }{\it ( blow-up)}: 

\smallskip

Il existe une suite de points $ (y_i)_i $, $ y_i \to x_0 $ et deux suites de r\'eels positifs $ (l_i)_i, (L_i)_i $, $ l_i \to 0 $, $ L_i \to +\infty $, telles qu'en posant $ v_i(z)=\dfrac{[y_i+z/[u_i(y_i)]^{2/(n-2)}]}{u_i(y_i)} $, on ait:

$$ 0 < v_i(z) \leq  \beta_i \leq 2^{(n-2)/2}, \,\, \beta_i \to 1. $$

$$  v_i(z)  \to \left ( \dfrac{1}{1+{|z|^2}} \right )^{(n-2)/2}, \,\, {\rm la convergence \,\, est \,\, uniforme \,\, sur \,\, tout \,\, compact \,\, de } \,\, {\mathbb R}^n . $$

$$ l_i^{(n-2)/2} u_i(y_i) \times \inf_{\Omega} u_i \to +\infty $$

\underbar {\bf Preuve de la proposition:}

\bigskip

On utilise $ (H) $, on peut supposer qu'il existe une suite  $ R_i>0, R_i \to 0 $ et $ c_i \to +\infty $, telles que,

$$ {R_i}^{(n-2)}(\sup_{B(x_0,R_i)} u_i)^{s} \inf_{\Omega} u_i \geq c_i \to +\infty, $$

Soit, $ x_i \in  { B(x_0,R_i)} $, tel que $ \sup_{B(x_0,R_i)} u_i=u_i(x_i) $ et $ s_i(x)=[R_i-|x-x_i|]^{(n-2)/2} u_i(x), x\in B(x_i, R_i) $. Alors, $ x_i \to x_0 $.

\bigskip

On a, 

$$ \max_{B(x_i,R_i)} s_i(x)=s_i(y_i) \geq s_i(x_i)={R_i}^{(n-2)/2} u_i(x_i)\geq \sqrt {c_i}  \to + \infty. $$ 

\bigskip

On pose :

$$ l_i=R_i-|x-x_i|,\,\, \bar u_i(y)= u_i(y_i+y),\,\,  v_i(z)=\dfrac{u_i \left (y_i+ z/[u_i(y_i)]^{2/(n-2)} \right ) } {u_i(y_i)}. $$

Il est clair que, $ y_i \to x_0 $. On obtient aussi:

$$ L_i= \dfrac{l_i}{(c_i)^{1/2(n-2)}} [u_i(y_i)]^{2/(n-2)}=\dfrac{[s_i(y_i)]^{2/(n-2)}}{c_i^{1/2(n-2)}}\geq \dfrac{c_i^{1/(n-2)}}{c_i^{1/2(n-2)}}=c_i^{1/2(n-2)}\to +\infty. $$

\bigskip

Si $ |z|\leq L_i $, alors $ y=[y_i+z/ [u_i(y_i)]^{2/(n-2)}] \in B(0,\delta_i l_i) $ avec $ \delta_i=\dfrac{1}{(c_i)^{1/2(n-2)}} $ et $ |y-y_i| < R_i-|y_i-x_i| $, d'o\`u, $ |y-x_i| < R_i $ et donc, $ s_i(y)\leq s_i(y_i) $, ce qui revient \`a \'ecrire,

$$ u_i(y) [R_i-|y-y_i|]^{(n-2)/2} \leq u_i(y_i) (l_i)^{(n-2)/2}. $$

Comme, $ |y-y_i| \leq \delta_i l_i $, $ R_i >l_i$ et $ R_i-|y- y_i| \geq R_i-\delta_i l_i>l_i-\delta_i l_i=l_i(1-\delta_i) $, on obtient,

$$ 0 < v_i(z)=\dfrac{u_i(y)}{u_i(y_i)} \leq \left [ \dfrac{l_i}{l_i(1-\delta_i)} \right ]^{(n-2)/2}\leq 2^{(n-2)/2} . $$

On pose alors, $ \beta_i=\left ( \dfrac{1}{1-\delta_i} \right )^{(n-2)/2} $, il est clair que $ \beta_i \to 1 $.

\bigskip

La fonction $ v_i $ v\'erifie l'\'equation suivante:

$$ \Delta v_i=\tilde V_i v_i^{(n+2)/(n-2)}+\dfrac{\tilde W_i}{[u_i(y_i)]^{2/(n-2)}} v_i^{n/(n-2)} $$

avec,

$$ \tilde V_i(z)=V_i \left [y_i+\dfrac{z}{[u_i(y_i)]^{2/(n-2)}} \right ] \,\,\, {\rm et} \,\,\, \tilde W_i(z)=W_i \left [y_i+\dfrac{z}{[u_i(y_i)]^{2/(n-2)}} \right ]. $$

En, utilisant les estimations elliptiques, les th\'eor\`emes d'Ascoli et de Ladyzenskaya, $ ( v_i)_i $ converge uniform\'ement sur tout compact  vers une fonction $ v $ solution sur $ {\mathbb R}^n $ de, 

$$ \Delta v = V(0) v^{N-1}, \,\, v(0)=1,\,\, 0 \leq v\leq 1\leq 2^{(n-2)/2}, $$

Sans nuire \`a la g\'en\'eralit\'e, on peut supposer que $ V(0)=n(n-2) $. 

\bigskip

Par le principe du maximum, on a $ v>0 $ sur $ {\mathbb R}^n $ et un r\'esultat de Caffarelli-Gidas-Spruck ( voir [C-G-S]) donne, $ v(z)=\left ( \dfrac{1}{1+{|z|^2}} \right )^{(n-2)/2} $. On obtient les m\^emes propri\'et\'es de convergence des $ v_i $ que dans un article pr\'ec\'edent (voir [B ]). La propostion 2 est prouv\'ee.

\bigskip

\underbar {\bf Coordonn\'ees Polaires} {\it (M\'ethode "Moving-Plane")}

\bigskip

Posons pour $ t\in ]-\infty, \log 2 ] $ et $ \theta \in {\mathbb S}_{n-1}
$ :

\bigskip

$ w_i(t,\theta)=e^{(n-2)t/2}u_i(y_i+e^t\theta), \,\,
\bar V_i(t,\theta)=V_i(y_i+e^t\theta) \,\,\, {\rm et }\,\,\,
\bar W_i(t,\theta)=W_i(y_i+e^t\theta). $

\bigskip

Par ailleurs, soit $ L $ l'op\'erateur $
L=\partial_{tt}-\Delta_{\sigma}-\dfrac{(n-2)^2}{4} $, avec $
\Delta_{\sigma} $ op\'erateur de Laplace-Baltrami sur $ {\mathbb
  S}_{n-1} $. 

\bigskip

La fonction $ w_i $  est solution de l'\'equation suivante :

$$ -Lw_i=\bar V_i{w_i}^{N-1}+e^t \times \bar W_i
{w_i}^{n/(n-2)}. $$

On pose pour $ \lambda \leq 0 $ :

\bigskip

$ t^{\lambda}=2\lambda-t  $ $ w_i^{\lambda}(t,\theta)=w_i(t^{\lambda},\theta) $, $ \bar
V_i^{\lambda}(t,\theta)=\bar V_i(t^{\lambda},\theta) $ et $ \bar
W_i^{\lambda}(t,\theta)=\bar W_i(t^{\lambda},\theta) .$

\bigskip

Alors, pour pouvoir v\'erifier si le Lemme 2 du Th\'eor\`eme 1 dans [B] reste
valable, il suffit de noter que  la quantit\'e $ -L( w_i^{\lambda}-w_i) $
est n\'egative lorsque $ w_i^{\lambda}-w_i $ l'est. En fait, pour chaque
indice $ i $, $ \lambda=
\xi_i \leq \log \eta_i+2 $, ($ \eta_i=[u_i(y_i)]^{(-2)/(n-2)}) $. 

\bigskip

Tout d'abord:

$$ w_i(2\xi_i-t,\theta)=w_i[(\xi_i-t+\xi_i-\log\eta_i-2)+(\log
\eta_i+2)] , $$

par d\'efinition de $ w_i $ et pour $ \xi_i \leq t $:

 $$ w_i(2\xi_i-t,\theta)=e^{[(n-2)(\xi_i-t+\xi_i-\log\eta_i-2)]/2}e^{n-2}v_i[\theta e^2e^{(\xi_i-t)+(\xi_i-\log\eta_i-2)}]
\leq 2^{(n-2)/2}e^{n-2}=\bar c. $$

On sait que

$$ -L( w_i^{\xi_i}-w_i)=[\bar V_i^{\xi_i
  }(w_i^{\xi_i})^{N-1}-\bar V_i
{w_i}^{N-1}]+[e^{t^{\xi_i}}{\bar W_i}^{\xi_i
  }(w_i^{\xi_i})^{n/(n-2) }-e^t\bar W_i
{w_i}^{n/(n-2)} ] ,  $$

Les deux termes du second membre, not\'es  $ Z_1 $ et $ Z_2 $, peuvent s'\'ecrire:

$$ Z_1=(\bar V_i^{\xi_i }-\bar V_i)(w_i^{\xi_i })^{N-1}+\bar
V_i[(w_i^{\xi_i })^{N-1}-{w_i}^{N-1}],$$

et

$$ Z_2=(\bar
W_i^{\xi_i }-\bar W_i)(w_i^{\xi_i })^{n/(n-2) }e^{
  t^{\xi_i }}+e^{
  t^{\xi_i }}\bar W_i[(w_i^{\xi_i })^{n/(n-2) }-{w_i}^{n/(n-2)
  }]+\bar W_i {w_i}^{n/(n-2) }(e^{t^{\xi_i }}-e^t ). $$

D'autre part, comme dans la d\'emonstration du Th\'eor\`eme 2 dans [B]:

$$ {w_i}^{\xi_i} \leq w_i \,\,\,{\rm et } \,\,\,
w_i^{\xi_i}(t,\theta)\leq \bar c \,\,\, {\rm pour \, tout } \,\,\,
(t,\theta)\in [\xi_i,\log 2] \times {\mathbb S}_{n-1}  , $$

o\`u  $ \bar c $ est une constante positive ind\'ependante de $ i
$ de $ w_i^{\xi_i} $ pour $ \xi_i \leq \log \eta_i+2 $;

$$ |\bar V_i^{\xi_i }-\bar V_i|\leq A_i (e^t-e^{ t^{\xi_i }}) \,\,\,
{\rm et } \,\,\, |\bar W_i^{\xi_i }-\bar W_i|\leq B_i (e^t-e^{ t^{\xi_i
    }}), $$

D'o\`u

\bigskip

$ Z_1 \leq A_i\, ({w_i^{\xi_i}})^{N-1} \,  (e^t-e^{ t^{\xi_i }}) \,\,\,
{\rm et} \,\,\,  Z_2 \leq  B_i \, ({(w_i ^{\xi_i})}^{n/(n-2) }\,  (e^t-e^{
  t^{\xi_i }})+ c\, {(w_i^{\xi_i})}^{n/(n-2)} \times   (e^{t^{\xi_i
    }}-e^t ) 
  $.

\bigskip

Ainsi, 

\bigskip

$ -L(w_i^{\xi_i}-w_i) \leq (w_i^{\xi_i})^{n/(n-2)}[ (A_i\,  {w_i^{\xi_i}}^{2/(n-2)}+ B_i)  \,  (e^t-e^{ t^{\xi_i }})+ c\,\,  (e^{ t^{\xi_i
    }}-e^t)].
$

\bigskip

Puisque $ w_i^{\xi_i} \leq \bar c $, on obtient:

\bigskip

$ -L(w_i^{\xi_i}-w_i)\leq  (w_i^{\xi_i})^{n/(n-2) } [(A_i {\bar c}
^{2/(n-2)}+ B_i)  \,  (e^t-e^{ t^{\xi_i }})+ c\,\,  (e^{t^{\xi_i
    }}-e^t ) ]. \,\,\,(1)$

\bigskip

D\'eterminons le signe de 
 $ \bar Z=[( A_i{\bar c}
^{2/(n-2)}+ B_i)  \, (e^t-e^{ t^{\xi_i }})+ c\,\,  (e^{ t^{\xi_i
    }}-e^t ) ]. $

\bigskip

L'in\'egalit\'e $ (1) $ devient alors :

$$ -L(w_i^{\xi_i}-w_i) \leq (w_i^{\xi_i})^{\alpha}[-c+ A_i \,{\bar c}^{2/(n-2)}+B_i]( e^t-e^{ t^{\xi_i }}). $$

On sait que $ A_i \to 0 $ et $ B_i \to 0 $. Pour $ t_0 < 0 $, assez petit, la quantit\'e $ c - A_i \,{\bar c}^{2/(n-2)}-B_i  $ devient
  positive et le r\'esultat cherch\'e est obtenu dans l'intervalle $
  [\xi_i,t_0] $.

\bigskip

Le fait de prendre l'intervalle  $  [\xi_i,t_0] $ au lieu de  $
[\xi_i, \log 2] $, n'est pas g\^enant, au contraire, plus l'intervalle
est petit plus l'infimum est grand. La suite de la d\'emonstration est
identique \'a celle de la fin du Th\'eor\`eme 1.

\bigskip

On pourrait croire que $ t_0 $ d\'epend de $ \xi_i $ ou de $
w_i^{\xi_i} $, mais  $ t_0 $ d\'epend seulement de $
\bar c $, une constante qui ne d\'epend que de $ n $, $ a $ et $ b $.

\bigskip

On calcule $ t_0 $ puis on introduit $ \xi_i \leq \log \eta_i+2 $ comme dans les autres
th\'eor\`emes, et on v\'erifie  l'in\'egalit\'e $ L(w_i^{\xi_i}-w_i )\leq 0 $, d\`es
que $ w_i^{\xi_i}-w_i \leq 0 $ sur $ [\xi_i,t_0].$

\bigskip

Ayant d\'etermin\'e $ t_0 <0 $ tel que $ c- A_i \,{\bar c}^{N-1-\alpha}-B_i $ soit positive, on pose:

\bigskip

$ \xi_i =\sup \{ \mu_i \leq \log \eta_i+2, w_i
  ^{\mu_i}(t,\theta)-w_i(t,\theta)  \leq 0, \forall \, (t,\theta) \in
  [\mu_i,t_0]\times {\mathbb S}_{n-1} \}$.

\bigskip

Par d\'efinition de $ \xi_i $, $ w_i^{\xi_i}-w_i \leq 0 $. Ensuite, on
 v\'erifie que $ -L(w_i ^{\xi_i}-w_i) \leq 0 $.

\bigskip

Comme dans le Th\'eor\`eme 1 dans [B], le principe du maximum, entra\^ine:

\bigskip

$ \qquad \qquad \min_{\theta \in {\mathbb S}_{n-1}}w_i(t_0,\theta) \leq \max_{\theta
  \in {\mathbb S}_{n-1}} w_i(2\xi_i-t_0) \,.$

\bigskip

Or,

$ \qquad w_i(t_0,\theta)=e^{t_0} u_i(a_i+e^{t_0}\theta)\geq e^{t_0} \min u_i
\,\,{\rm et} \,\, w_i(2\xi_i-t_0)\leq \dfrac{ c_0 }{u_i(a_i)} $,

\bigskip

donc:

$$ u_i(a_i) \times \min u_i \leq c .$$

Ce qui contredit la proposition.

\bigskip

\underbar {\bf Preuve du Th\'eor\`eme 2}

\bigskip

Les \'etapes sont identiques \`a celles de la preuve du th\'eor\`eme 1. Il y a quelques modifications dans la partie " Coordonn\'ees polaires et m\'ethode moving-plane". La proposition de la preuve du th\'eor\`eme 1 se conserve. Notons que la technique blow-up se simplifie car $ u_i $ est d\'ecroissante et son maximum est atteint en 0.

\bigskip

\underbar {\bf Coordonn\'ees polaires} {\it (M\'ethode "Moving-plane")}

\bigskip

\bigskip

Posons pour $ t\in ]-\infty, \log 2 ] $ et $ \theta \in {\mathbb S}_{n-1}
$ :

\bigskip

$ w_i(t,\theta)=e^{(n-2)t/2}u_i(e^t), \,\,
\bar V_i(t,\theta)=V_i(e^t) \,\,\, {\rm et }\,\,\,
\bar W_i(t,\theta)=W_i(e^t). $

\bigskip

Par ailleurs, soit $ L $ l'op\'erateur $
L=\partial_{tt}-\dfrac{(n-2)^2}{4} $.

\bigskip

La fonction $ w_i $  est solution de l'\'equation suivante :

$$ -Lw_i=\bar V_i{w_i}^{N-1}+e^{2t} \bar W_i w_i. $$

On pose pour $ \lambda \leq 0 $ :

\bigskip

$ t^{\lambda}=2\lambda-t  $, $ w_i^{\lambda}(t,\theta)=w_i(t^{\lambda}) $, $ \bar
V_i^{\lambda}(t,\theta)=\bar V_i(t^{\lambda}) $ et $ \bar
W_i^{\lambda}(t,\theta)=\bar W_i(t^{\lambda}) .$

\bigskip

Alors, pour pouvoir v\'erifier si le Lemme 2 du Th\'eor\`eme 1 dans [B] reste
valable, il suffit de noter que  la quantit\'e $ -L( w_i^{\lambda}-w_i) $
est n\'egative lorsque $ w_i^{\lambda}-w_i $ l'est. En fait, pour chaque
indice $ i $, $ \lambda=
\xi_i \leq \log \eta_i+2 $, ($ \eta_i=[u_i(y_i)]^{(-2)/(n-2)}) $. 

\bigskip

Tout d'abord:

$$ w_i(2\xi_i-t)=w_i[(\xi_i-t+\xi_i-\log\eta_i-2)+(\log
\eta_i+2)] , $$

par d\'efinition de $ w_i $ et pour $ \xi_i \leq t $:

 $$ w_i(2\xi_i-t)=e^{[(n-2)(\xi_i-t+\xi_i-\log\eta_i-2)]/2}e^{n-2}v_i[e^2e^{(\xi_i-t)+(\xi_i-\log\eta_i-2)}]
\leq 2^{(n-2)/2}e^{n-2}=\bar c. $$

On sait que

$$ -L( w_i^{\xi_i}-w_i)=[\bar V_i^{\xi_i
  }(w_i^{\xi_i})^{N-1}-\bar V_i
{w_i}^{N-1}]+[e^{2t^{\xi_i}}{\bar W_i}^{\xi_i
  }(w_i^{\xi_i})-e^{2t}\bar W_i w_i ] ,  $$

Les deux termes du second membre, not\'es  $ Z_1 $ et $ Z_2 $, peuvent s'\'ecrire:

$$ Z_1=(\bar V_i^{\xi_i }-\bar V_i)(w_i^{\xi_i })^{N-1}+\bar
V_i[(w_i^{\xi_i })^{N-1}-{w_i}^{N-1}],$$

et

$$ Z_2=[(e^{2t}\bar W_i)^{\xi_i}-(e^{2t}\bar W_i)]w_i^{\xi_i} + e^{2t}\bar W_i(w_i^{\xi_i } -{w_i}). $$

D'autre part, comme dans la d\'emonstration du Th\'eor\`eme 2 dans [B]:

$$ {w_i}^{\xi_i} \leq w_i \,\,\,{\rm et } \,\,\,
w_i^{\xi_i}(t,\theta)\leq \bar c \,\,\, {\rm pour \, tout } \,\,\,
(t,\theta)\in [\xi_i,\log 2] \times {\mathbb S}_{n-1}  , $$

o\`u  $ \bar c $ est une constante positive ind\'ependante de $ i
$ de $ w_i^{\xi_i} $ pour $ \xi_i \leq \log \eta_i+2 $;

$$ |\bar V_i^{\xi_i }-\bar V_i|\leq A_i (e^{2t}-e^{2t^{\xi_i }}) \,\,\,
{\rm et } \,\,\, |(e^{2t}\bar W_i)^{\xi_i }-(e^{2t}\bar W_i)-W_i(0)(e^{2t^{\xi_i}}-e^{2t})|\leq \tilde B_i (e^{2t}-e^{ 2t^{\xi_i
    }}), $$

avec, $ \tilde B_i \to 0 $. D'o\`u

\bigskip

$ Z_1 \leq A_i\, ({w_i^{\xi_i}})^{N-1} \,  (e^{2t}-e^{2t^{\xi_i }}) \,\,\,
{\rm et} \,\,\,  Z_2 \leq  \tilde B_i \, (w_i^{\xi_i})\,  (e^{2t}-e^{2t^{\xi_i }})+ c\, (w_i^{\xi_i}) \times   (e^{2t^{\xi_i
    }}-e^{2t}) 
  $.

\bigskip

Ainsi, 

\bigskip

$ -L(w_i^{\xi_i}-w_i) \leq w_i^{\xi_i} [ [A_i\,  ({w_i^{\xi_i}})^{4/(n-2)}+ \tilde B_i]  \,  (e^{2t}-e^{2t^{\xi_i }})+ c\,\,  (e^{2t^{\xi_i}}-e^{2t})].
$

\bigskip

Puisque $ w_i^{\xi_i} \leq \bar c $, on obtient:

\bigskip

$ -L(w_i^{\xi_i}-w_i)\leq  w_i^{\xi_i} [(A_i {\bar c}
^{4/(n-2)}+ \tilde B_i)  \,  (e^{2t}-e^{2t^{\xi_i }})+ c\,\,  (e^{2t^{\xi_i}}-e^{2t} ) ]. \,\,\,(1)$

\bigskip

D\'eterminons le signe de 
 $ \bar Z=[( A_i{\bar c}
^{4/(n-2)}+ \tilde B_i)  \, (e^{2t}-e^{2t^{\xi_i }})+ c\,\,  (e^{2t^{\xi_i}}-e^{2t}) ]. $

\bigskip

L'in\'egalit\'e $ (1) $ devient alors :

$$ -L(w_i^{\xi_i}-w_i) \leq w_i^{\xi_i}[-c+ A_i \,{\bar c}^{4/(n-2)}+\tilde B_i]( e^{2t}-e^{2t^{\xi_i }}). $$

On sait que $ A_i \to 0 $ et $ \tilde B_i \to 0 $. Pour $ t_0 < 0 $, assez petit, la quantit\'e $ c - A_i \,{\bar c}^{4/(n-2)}-\tilde B_i  $ devient positive et le r\'esultat cherch\'e est obtenu dans l'intervalle $
  [\xi_i,t_0] $.

\bigskip

Le fait de prendre l'intervalle  $  [\xi_i,t_0] $ au lieu de  $
[\xi_i, \log 2] $, n'est pas g\^enant, au contraire, plus l'intervalle
est petit plus l'infimum est grand. La suite de la d\'emonstration est
identique \'a celle de la fin du Th\'eor\`eme 1.

\bigskip

On pourrait croire que $ t_0 $ d\'epend de $ \xi_i $ ou de $
w_i^{\xi_i} $, mais  $ t_0 $ d\'epend seulement de $
\bar c $, une constante qui ne d\'epend que de $ n $, $ a $ et $ b $.

\bigskip

On calcule $ t_0 $ puis on introduit $ \xi_i \leq \log \eta_i+2 $ comme dans les autres
th\'eor\`emes, et on v\'erifie  l'in\'egalit\'e $ L(w_i^{\xi_i}-w_i )\leq 0 $, d\`es
que $ w_i^{\xi_i}-w_i \leq 0 $ sur $ [\xi_i,t_0].$

\bigskip

Ayant d\'etermin\'e $ t_0 <0 $ tel que $ c- A_i \,{\bar c}^{4/(n-2)}-\tilde B_i $ soit positive, on pose:

\bigskip

$ \xi_i =\sup \{ \mu_i \leq \log \eta_i+2, w_i
  ^{\mu_i}(t)-w_i(t)  \leq 0, \forall \, t \in
  [\mu_i,t_0] \}$.

\bigskip

Par d\'efinition de $ \xi_i $, $ w_i^{\xi_i}-w_i \leq 0 $. Ensuite, on
 v\'erifie que $ -L(w_i ^{\xi_i}-w_i) \leq 0 $.

\bigskip

Comme dans le Th\'eor\`eme 1 dans [B], le principe du maximum, entra\^ine:

$$ w_i(t_0) \leq  w_i(2\xi_i-t_0), $$

comme $ u_i $ est d\'ecroissante, on obtient:

$$ u_i(a_i) \times u_i(1) \leq c .$$

\bigskip

\underbar {\bf R\'ef\'erences:}

\bigskip

[A] T. Aubin. Nonlinear Problems in Riemannian Geometry. Springer-Verlag 1998.

\smallskip

[B] S.S Bahoura. Majorations du type $ \sup u \times \inf u \leq c $ pour l'\'equation de la courbure scalaire prescrite sur un ouvert de $ {\mathbb R}^n, n\geq 3 $. J.Math.Pures Appl.(9) 83 (2004), no.9, 1109-1150.

\smallskip

[C-G-S] Caffarelli L, Gidas B., Spruck J. Asymptotic symmetry and local
behavior of semilinear elliptic equations with critical Sobolev
growth. Commun. Pure Appl. Math. 37 (1984) 369-402.

\smallskip

[C-L 1] C.C. Chen, C-S. Lin. Prescribing scalar curvature on $ {\mathbb S}_n $. I. A priori estimates.  J. Differential Geom.  57  (2001),  no. 1, 67--171.

\smallskip

[C-L 2] Chen C-C. and Lin C-S. Blowing up with infinite energy of conformal
metrics on $ {\mathbb S}_n $. Comm. Partial Differ Equations. 24 (5,6)
(1999) 785-799.

\smallskip

[G-N-N] B. Gidas, W. Ni, L. Nirenberg, Symmetry and Related Propreties via the Maximum Principle, Comm. Math. Phys., vol 68, 1979, pp. 209-243.

\smallskip

[L1] Y.Y Li. Prescribing Scalar Curvature on $ {\mathbb S}_n $ and related Problems. I. J. Differential Equations 120 (1995), no. 2, 319-410.

\smallskip

[L2] Y.Y Li. Prescribing Scalar Curvature on $ {\mathbb S}_n $ and related Problems. II. Comm. Pure. Appl. Math. 49(1996), no.6, 541-597.

\end{document}